\documentclass[10pt,notitlepage,twoside,a4paper]{amsart}
 \usepackage{amsfonts}
\usepackage[latin1]{inputenc}
\usepackage[cyr]{aeguill}
\usepackage{xspace}
\usepackage[francais]{babel}

\usepackage{amsmath,amssymb,enumerate}

\usepackage{epsfig,fancyhdr,color}
 
\usepackage[table]{xcolor}
 
\usepackage{amssymb}
\usepackage{amsmath,amsthm}
\usepackage{latexsym}
\usepackage{amscd}
\usepackage{psfrag}
\usepackage{graphicx}
\usepackage[latin1]{inputenc}
\usepackage[all]{xy}
\usepackage[mathcal]{eucal}

\definecolor{NoteColor}{rgb}{1,0,0}

\renewcommand{\textsc}{\textcolor{red}}

\newtheorem*{theorem 1}{\rm\bf Proposition 1}
\newtheorem*{theorem 2}{\rm\bf Proposition 2}

\theoremstyle{definition}

\theoremstyle{remark}

\def\interieur#1{\mathord{\mathop{\kern 0pt #1}\limits^\circ}}

\title[Messiaen et les mathématiques]{Messiaen et les mathématiques}

\author[Athanase Papadopoulos]{Athanase Papadopoulos}

 \address{Athanase Papadopoulos,  Institut de Recherche Math\'ematique Avanc\'ee\\ CNRS et Universit\'e de Strasbourg\\\small 7 rue Ren\'e
  Descartes - 67084 Strasbourg Cedex, France et Laboratoire GREAM (Groupe de recherches expérimentales sur l'acte musical), 
Université de Strabsourg,  5 rue du Général Rouvillois - 67083 Strasbourg Cedex}

\date{}

\begin{document}
\maketitle
 
\begin{abstract}
Nous passons en revue l'utilisation de certaines notions mathématique dans les compositions d'Olivier Messiaen et dans son travail théorique

La version finale de cet article a paru dans l'ouvrage 
\emph{Twentieth-Century Music and Mathematics},
R. Illiano (ed.), Brepols, Turnhout, 2019.

\end{abstract}

\medskip

\medskip

\begin{quote}\smaller
Abstract. We review some of Olivier Messiaen's use of mathematics in his composition and his theoretical writings.

The final version of this paper appeared in the book \emph{Twentieth-Century Music and Mathematics},
R. Illiano (ed.), Brepols, Turnhout, 2019.

\medskip

\medskip 

Mots clés :   mathematiques et musique, rythme, modes à transpositions limitées, rythme non rétrogradable, contrepoint, permutations symétrique, Olivier Messiaen, 

Classification AMS: 00A65 (Mathematics and music)

\larger
\end{quote}

\medskip 

\medskip 

\larger
\section{Introduction}

Olivier Messiaen a marqué le vingtième siècle par un langage rythmique et modal original et unique, qu'il a transmis à travers ses compositions, ses écrits et ses cours au Conservatoire de Paris. Son \oe uvre, hautement novatrice, est en même temps ancrée dans des traditions très anciennes, tout particulièrement celles de la musique de l'antiquité grecque, du plein-chant médiéval et de la musique hindoue. 

Notre sujet principal dans cet article concerne les mathématiques dans l'\oe uvre de Messiaen, et il est peut-être utile de rappeler ici qu'en mathématiques, la découverte est toujours ancrée dans une longue tradition, souvent bi-millénaire. Il n'est donc pas étonnant que cela soit aussi le cas en musique, un domaine partageant avec les mathématiques une histoire et un mode de pensée communes, et cela d'autant plus pour la musique de Messiaen, qui a un caractère mathématique prononcé. 

Messiaen est né  en 1908, c'est-à-dire huit années après le début du dernier siècle, et il est mort en 1992,  huit années avant sa fin ;  le quadruplet 1900, 1908, 1992, 2000 donne lieu au rythme 
(8, 84, 8) qui possède une symétrie qui en fait (ce que Messiaen appelle) un
rythme non rétrogradable, une notion dont on parlera plus bas.

Sans prétendre rendre compte tant soit peu de l'importance de l'\oe uvre gigantesque de Messiaen, je me suis limité sans cet article à certains aspects math\'ematiques de cette \oe uvre. Ils sont tous basés sur des notions élémentaires : nombres premiers,
permutations et autres sym\'etries s'exprimant naturellement en termes de théorie de groupes (même si l'on ne sait pas si Messiaen a explicitement utilisé ce vocabulaire, synonyme de symétrie).  
Ces notions, bien que très simples, sont au c\oe ur des mathématiques -- les mathématiciens professionnels ont horreur de la complication -- et elles ont l'avantage d'être intemporelles, c'est-à-dire en dehors des goûts du jour. Il ne s'agit, par exemple, ni de musique fractale, ni de musique quantique, ni de musique aléatoire.
 M\^eme si
 Messiaen ne s'est jamais consid\'er\'e comme un math\'ematicien -- et, à ce propos, il est légitime de se poser la question : qu'est-ce qu'un mathématicien ? -- il a utilisé dans sa musique les nombres, les
suites de  nombres, leurs transformations, et leurs sym\'etries de façon consciente et systématique et il les a mis en valeur dans ses écrits. Plus généralement, il a 
 accord\'e aux mathématiques
une place de premier ordre dans ses compositions et dans son \oe uvre
th\'eorique.    

Au-del\`a de la discussion sur Messiaen d\'evelopp\'ee
dans cet article, l'un des th\`emes que je voudrais illustrer 
est que la musique contribue à donner vie et à rendre perceptibles \`a nos sens des notions dont la nature serait \emph{a priori} purement abstraite. Enfin, l'article sera aussi l'occasion de passer en revue, même si c'est de façon très brève,  certains côtés historiques, philosophiques et théologiques de la pensée de Messiaen.

\section{Rythme}

 Le rythme occupe une place de premier plan dans la musique de Messiaen. 
   Dans un entretien avec le journaliste et critique musical Claude Samuel, publié en 1999, il déclare : 
 « Je consid\`ere que le rythme est la
partie primordiale et peut-\^etre essentielle de la musique ; je pense qu'il a
vraisemblablement exist\'e avant la m\'elodie et l'harmonie, et j'ai enfin une
pr\'ef\'erence secr\`ete pour cet \'el\'ement »\footnote{Samuel 1999, p. 101.}. L'\oe uvre
   th\'eorique monumentale, sur laquelle Messiaen travailla pendant 
40 ans et qui fut publiée de manière posthume, s'intitule « Traité de rythme, de couleur, et d'ornithologie ».  Dans ce titre, le rythme vient en premier lieu, et le premier tome du traité (375 pages), est entièrement dédié à cette notion. 
 Dans les premi\`eres
pages de ce volume, Messiaen s'\'el\`eve contre l'opinion commune qui dit que
la musique est faite avec des sons :  « Je dis non ! Non, pas seulement avec des
sons... La m\'elodie n'existerait pas sans le Rythme!... La musique est faite
d'abord avec des {\it Dur\'ees}, des {\it \'Elans} et des {\it Repos}, des
{\it Accents}, des {\it Intensit\'es} et des {\it Densit\'es}, des {\it
Attaques} et des {\it Timbres}, toutes choses qui se regroupent sous un vocable
g\'en\'eral : le {\it  Rythme} »\footnote{Messiaen 2001, t. 1, p. 40}.

S'agissant de rythme, Messiaen ne pensait pas à des rythmes réguliers, mais plutôt aux rythmes libres de la Grèce antique et de l'Inde, à ceux des vagues de la mer, des cascades en montagne, et du vent. Dans un de ses entretiens avec Claude Samuel, il dit: « Je ne fais aucune limitation entre le bruit et le son ; tout cela représente toujours pour moi de la musique »\footnote{Samuel 1967, p. 28.}. Il est intéressant de noter que l'observation de tels phénomènes naturels a aussi inspiré l'un des plus grands mathématiciens du vingtième siècle, René Thom. Ce dernier raconte dans ses écrits\footnote{Thom 1988.} que c'est en regardant le déploiement d'une vague qu'il a eu l'idée de celui d'une courbe algébrique dans le plan, qui fut le ferment ce qui est devenu plus tard la théorie des catastrophes. Il se souvenait toujours du moment où il a eu pour la première fois cette idée : étendu sur sa couchette, pendant sa première traversée de l'Atlantique. On n'a pas suffisamment souligné la similarité des modes de pensée de Thom et de  Messiaen. Pour l'un comme pour l'autre de ces deux grands penseurs, l'observation des phénomènes du monde réel a servi de ferment aux découvertes les plus abstraites. Signalons au passage que tous deux partageaient une grande sensibilité hellénistique\footnote{Papadopoulos 2019A et 2019B.}.

 Les couleurs, deuxième notion à laquelle fait référence le titre du traité de Messiaen,  sont celles qu'il associait aux sons et aux combinaisons de sons. Un exemple que l'on donne souvent de musique colorée est celui de la \emph{Cité Céleste}, une \oe uvre pour piano, ensemble à vent et percussions qu'il composa en 1963, dont le titre fait référence au livre de l'Apocalypse et qui, selon Messiaen, représente « une rosace de cathédrale aux couleurs flamboyantes et invisibles »\footnote{Note d'Olivier Messiaen dans Messiaen 1966.}.
 
 Enfin, l'ornithologie, dans le titre du traité,  fait référence aux oiseaux que Messiaen passait son temps à observer et dont il notait les chants : encore un matériau fondamental de la nature avec lequel il travaillait. Messiaen dit à plusieurs occasions qu'il considérait ces animaux comme les musiciens les plus doués qui existent sur
notre plan\`ete, et surtout, les meilleurs rythmiciens\footnote{Messiaen 2001, t. 1, p. 53.}. Ils étaient  ses maîtres en musique et il leur rendit hommage dans de multiples compositions. Dans son dialogue avec Samuel, il déclare : « Vous ne trouverez jamais, dans le chant d'oiseaux, une erreur
de rythme, de m\'elodie ou de contrepoint »\footnote{Samuel 1999.}. Les oiseaux, disait-il souvent, ont tout
invent\'e, y compris l'improvisation collective. « Tout ce que je sais de la m\'elodie, ce sont les
oiseaux qui me l'ont appris...  »\footnote{Messiaen 2001, t. 1, p. 53.}. Et sur le rythme : « La Grive
musicienne -- qui est peut-\^etre la reine des oiseaux chanteurs -- poss\`ede
un chant magique, incantatoire, coup\'e en petites formules rythmiques tr\`es
nettes, toujours r\'ep\'et\'ees de 2 \`a 5 fois, le plus souvent 3 fois (comme
dans le rituel des invocations religieuses et les enchantements de la
sorcellerie primitive). En dehors de quelques rythmes caract\'eristiques, les
strophes sont toujours nouvelles et l'invention rythmique in\'epuisable.
L'agencement des dur\'ees et des nombres, toujours inattendu, impr\'evu,
surprenant, manifeste cependant un sens de l'\'equilibre tel qu'on a peine \`a
croire \`a une improvisation. »\footnote{Messiaen 2001, t. 1, p. 54.}. L'analyse des \emph{Couleurs de la cité Céleste}, une \oe uvre que nous avons mentionnée,   révèle qu'elle fait intervenir les chants de 21 oiseaux exotiques\footnote{Mathon 1995.}.

Le rythme n'est pas \'etranger aux math\'ematiques, ne
serait-ce que parce que les dur\'ees, les repos, les intensit\'es, les densit\'es, etc. qui font partie des caractéristiques d'un rythme 
  se mesurent avec des nombres. 
  
  Pour aborder de manière plus précise certains aspects mathématiques du rythme, nous allons considérer ce dernier comme une suite de dur\'ees, c'est-\`a-dire comme un 
« d\'ecoupage du temps », expression qui prend toute sa valeur sous la plume de 
Messiaen, qui d'ailleurs commence son \emph{Trait\'e} par une longue
\'etude sur le Temps : temps absolu, temps relatif, temps biologique, temps
cosmologique, temps physiologique, temps psychologique, etc., avec de longues
digressions sur la notion de Temps dans la mythologie, dans la Bible, chez
Einstein, Bergson  et d'autres. Signalons en passant que Messiaen a utilisé la notion de rythme dans un sens beaucoup plus élargi, 
 qui n'est pas confiné à la musique ; il voyait le rythme dans la nature, pas seulement dans les sons mais aussi dans les formes.  Il le lisait par exemple dans les ailes des papillons.  Il écrit dans son \emph{Trait\'e de rythme, de couleur et
d'ornithologie} : « Il doit exister quelque part un musée des formes rythmiques, des archétypes de la branche, de la feuille et de la fleur, qui seraient en même temps l'origine de toutes les formes et de tous les rythmes connus »\footnote{Messiaen 2001, t. 1, p. 56.}.  Les «  formes » dans la nature, c'était aussi l'objet d'étude et de pensée par excellence de Thom, que l'on a mentionné plus haut.

\medskip 
 Messiaen se définissait d'abord comme rythmicien. Il aimait les rythmes de la Gr\`ece antique et de l'Inde, à cause de leurs propriétés arithmétiques, mais aussi et surtout pour la grande liberté que leur utilisation lui offrait. On est loin des rythmes monotones et réguliers à trois temps, à quatre temps, etc., et même si les rythmes grecs et hindous sont structurés et classés, le nombre de possibilités qu'ils offrent à l'intérieur d'une même composition est infini ; on n'est pas limité par une règle uniforme \emph{a priori}. Ils sont à l'image des mathématiques qui, loin de confiner la pensée dans un cadre délimité, lui fournissent un espace de liberté infinie. Messiaen déclare, à propos du \emph{Vent de l'Esprit},  dernière pièce de la \emph{Messe de la Pentecôte}, qu'il « mélange la chose la plus vivante, la plus libre qui soit : un chant d'Alouette -- avec une combinaison rythmique de la plus extrême rigueur  »\footnote{Cité dans Ide 1999.}.

\medskip 

 Quelques rappels sur les rythmes grecs sont peut-être utiles ici.

\medskip 
  
On sait que dans la Gr\`ece antique, la musique accompagnait la po\'esie. Ainsi, le
rythme musical  suivait celui de la d\'eclamation.  Il y avait deux
sortes de 
 dur\'ees pour les notes, des dur\'ees {\it
longues}, correspondant aux syllabes longues, et des dur\'ees {\it br\`eves}, correspondant aux syllabes br\`eves. 
Les {\it mesures} (si l'on peut parler de mesure), y  sont de longueur in\'egale ; elles sont comparables aux respirations entre les phrases du langage parlé, ou déclamé.  C'est une libert\'e que le chant gr\'egorien avait préservée (« rythmique verbale »)  et  que Messiaen utilisa dans toutes ses compositions.
Le rythme en ce sens, comme succession de dur\'ees longues et
br\`eves, s'appelle   {\it m\`etre} ou {\it m\'etrique}. 
 Les diff\'erents types de m\'etriques grecques furent analysés et class\'es  par les théoriciens grecs de la musique, en particulier par  Aristox\`ene de Tarente, le grand
musicologue du 4\`eme si\`ecle av. J. C.  Messiaen est 
probablement le principal compositeur du 20e si\`ecle 
\`a avoir sérieusement raviv\'e l'usage d'une telle musique, dans ses
compositions ainsi que dans son enseignement 
th\'eorique, même si, comme il dit lui-m\^eme, cette tradition
n'avait pas \'et\'e compl\`etement perdue dans la musique occidentale, puisqu'on en 
retrouve des traces dans la
musique folklorique roumaine et chez certains compositeurs comme Ravel ou 
Stravinski.   Le tome 1 du \emph{Trait\'e de Rythme, de couleur et d'ornithologie} contient un chapitre de 170
pages sur la m\'etrique grecque. Par son usage du rythme, Messiaen se sentait proche des Orientaux : « Les Orientaux sont tous rythmiciens, les Hindous plus que tous les autres réunis. Les Occidentaux sont plus harmonistes que rythmiciens »\footnote{Messiaen 2001, t. 1, p. 30.}.

 \medskip 

Une autre  caract\'eristique de la m\'etrique grecque, qui
\'etait rarement utilisée dans  la musique classique occidentale quand Messiaen a entrepris da formation, est qu'elle utilise souvent des rythmes dont la dur\'ee totale est un nombre
premier (en particulier 5, 11 et 17) . Un
  exemple  de rythme \`a 5 temps est le rythme   {\it Cr\'etique} (ou {\it
Cr\'etois}), représenté par la suite 2, 1, 2 (c'est-\`a-dire deux longues, une
br\`eve et deux longues) et ses {\it permutations} : 2, 2, 1 et 1, 2, 2. Cet
exemple nous introduit directement  \`a deux notions importantes dans le  monde
musical  de Messiaen. La premi\`ere est manifest\'ee par  la
sym\'etrie centrale de la suite  2, 1, 2,  qui en fait 
  un exemple de ce que Messiaen appelle {\it rythme non r\'etrogradable}, que l'on a déjà mentionné dans l'introduction de cet article, et dont on
reparlera encore plus loin. Le rythme 2, 1, 2  s'appelle {\it Amphimacros}, qui
veut dire (comme l'explique  Messiaen) « longue entrourant la br\`eve ». L'autre 
notion, sur laquelle on reviendra aussi plus loin, est celle de {\it
permutation}, que l'on a appliqué ici à une suite de  dur\'ees.

Passons maintenant aux rythmes hindous, et commençons par nous poser la question :
Pourquoi Messiaen s'est-il intéressé à la musique hindoue ?
 Il nous en donne lui-même des éléments de réponse. Il \'ecrit, dans le tome 1 de son \emph{Traité}\footnote{Messiaen 2001, t. 1, p. 258.} : « La musique
hindoue est celle qui est certainement all\'ee le plus loin dans le domaine
rythmique et sp\'ecialement dans l'ordre quantitatif (combinaisons des longues
et des br\`eves). Les rythmes hindous, d'un raffinement, d'une subtilit\'e sans
\'egals, laissent loin derri\`ere eux nos pauvres rythmes occidentaux avec leurs
mesures isochrones, et leurs perp\'etuelles divisions et multiplications par 2
(quelquefois par 3). »

\medskip 

Comme dans le cas des rythmes grecs, les rythmes hindous foisonnent dans 
les compositions de Messiaen, et leur étude occupe dans le tome 1 de son \emph{Traité} une
partie aussi importante que des  rythmes grecs. Messiaen y
dresse le tableau des 120 {\it De\c
c\^\i-T\^alas}\footnote{En hindou, De\c c\^\i = rythme ; T\^ala =
province. Ainsi, De\c c\^\i-T\^ala veut dire rythme des diff\'erentes provinces
(explication de Messiaen). Un catalogue de 120 talas avait été compilé par le musicologue indien du 13e siècle Sharngadeva ; c'est celui que Messiaen reprend. Ce dernier en avait déjà utilisé plusieurs, même avant de savoir qu'ils existaient dans la musique hindoue.}, des 36 rythmes de la radition karn\^atique, et
d'autres groupes de rythmes hindous.
Ces rythmes peuvent \^etre vus comme des suites de nombres ayant des propriétés arithmétiques, et ces dernières exerçaient une grande fascination sur Messiaen.  Dans ses entretiens avec
Pierre Samuel\footnote{Samuel 1999, p. 118.}, r\'epondant \`a une question sur les
  {\it De\c c\^\i -T\^alas}, Messiaen dit : « De toutes  fa\c cons, j'\'etais
orient\'e vers ces recherches, vers les divisions asym\'etriques, et vers un
\'el\'ement qui se rencontre dans la m\'etrique grecque et dans les rythmes de
l'Inde : les nombres premiers. Quand j'\'etais enfant, j'aimais d\'ej\`a les
nombres premiers, ces nombres qui, par le simple fait qu'ils ne sont pas
divisibles en fractions \'egales, d\'egagent une force occulte, puisque vous
savez que la Divinit\'e est indivisible... »
Le tableau des {\it De\c c\^\i-T\^alas} que Messiaen donne dans son \emph{Traité} foisonne de 
rythmes dont la somme des dur\'ees est un nombre premier (on trouve par exemple
les entiers  5, 7, 11, 17, 19, 37).  Il \'ecrit dans le tome 1, p.
266 : « L'impossibilit\'e de la division d'un nombre premier (autre que par
lui-m\^eme et par l'unit\'e) leur conf\`ere une sorte de puissance qui est
tr\`es effective dans le domaine du rythme. »

 \medskip
 
Une autre propri\'et\'e arithmétique, que l'on retrouve dans certains {\it De\c c\^\i-T\^alas} est celle
de 
 rythme consistant en une suite de dur\'ees suivie de son {\it augmentation},
par exemple la suite 1, 1, 1, 2, 2, 2 (T\^ala No. 73), qui est form\'e par la
suite 1, 1, 1, suivie de son augmentation par   multiplication par 2. Inversement, certains rythmes sont suivis par leur  {\it diminution}, parfois de nouveau par la  {\it diminution} de leur diminution, etc., comme dans l'exemple 4, 4, 2, 2, 1, 1 (T\^ala No. 115)

\medskip 

On trouve aussi des combinatoires plus complexes, par exemple le rythme  1, 3, 2,
3, 3, 3, 2, 3, 1, 3 (T\^ala No. 27) o\`u
 les dur\'ees d'ordre pair sont toutes \'egales entre elles, alors que les
dur\'ees d'ordre impair forment une suite croissante puis d\'ecroissante, de fa\c
con r\'eguli\`ere.

\medskip

On trouve aussi dans les  {\it De\c c\^\i-T\^alas} plusieurs exemples de  rythmes  non
r\'etrogradables, une propri\'et\'e que  l'on a rencontr\'ee dans la
m\'etrique grecque, c'est-\`a-dire des rythmes form\'es par une suite de dur\'ees
suivie de la suite sym\'etrique (avec \'eventuellement une valeur commune au
centre). On  retrouve
 par exemple le rythme grec {\it Amphimacre}, 2, 1, 2   (T\^ala No. 58), dont
Messiaen dit qu'il est « le plus simple et le plus naturel des rythmes non
r\'etrogradables, \'etant bas\'e sur le nombre 5, le nombre des doigts de la
main » (t. 1 p. 289).  D'autres   rythmes non r\'etrogradables  
dans les {\it De\c c\^\i-T\^alas} sont par exemple  2, 2, 1, 1, 2, 2 (T\^ala No.
26), ou 1, 1, 2, 2, 1, 1 (T\^ala No. 80), ou 2, 1, 1, 1, 1, 2 (T\^ala No. 111),
et il y en a plusieurs  d'autres. On reviendra plus loin sur l'importance de
ces rythmes.

  Messiaen a mis en valeur la notion de  rythme non r\'etrogradable tr\`es t\^ot dans
ses compositions, et   il y a attach\'e une importance capitale d\`es ses
premiers \'ecrits th\'eoriques (\emph{Technique de mon langage musical}, publi\'e en
1944)\footnote{Messiaen 1944.}. Un tel rythme est repr\'esent\'e par une suite de dur\'ees ayant la
propri\'et\'e que, quand on la lit de gauche \`a droite ou de droite \`a gauche,
on obtient la m\^eme suite. Il est bon de rappeler ici que la {\it
r\'etrogradation} est un  proc\'ed\'e  classique en   {\it  contrepoint}
(qui est l'art de
transformer et de superposer des lignes m\'elodiques). Ce proc\'ed\'e  consiste
\`a lire un certain   motif musical \`a l'envers (c'est-\`a-dire en commen\c
cant
par la derni\`ere note et en terminant par la premi\`ere).
 Le motif initial s'appelle alors {\it motif en mouvement direct}, et le transformé s'appelle {\it  motif r\'etrograde}. Ainsi, on
peut voir un rythme non
r\'etrogradable comme la juxtaposition d'un rythme en mouvement direct suivi de son
mouvement r\'etrograde (avec \'eventuellement une valeur commune au milieu). La
r\'etrogradation, comme technique de contrepoint,   est utilis\'ee et
enseign\'ee depuis le 14\`eme si\`ecle. Avant Messiaen, elle \'etait appliqu\'ee
\`a une ligne m\'elodique, c'est-à-dire une suite de hauteurs.  
Avec Messiaen, elle acquiert un caract\`ere plus abstrait ; elle est appliqu\'ee 
au rythme, ind\'ependamment de la hauteur des notes. L'auditeur de la musique de
Messiaen est appel\'e \`a sentir la r\'etrogradation d'un motif musical
uniquement au niveau des dur\'ees, sans qu'il y ait de
correspondance r\'eguli\`ere (égalité, transposition, sym\'etrie, etc.) entre  les
hauteurs des  notes (comme fr\'equences)
du motif en mouvement direct avec celles du motif en mouvement r\'etrograde.

\medskip 

Dans le tome 1 de son Trait\'e\footnote{Messiaen 2001, t. 1, p. 26.}, Messiaen, cite 8 mesures ssuccessives
de la « Danse de la fureur, pour les sept trompettes » (6e mouvement de son 
 {\it Quatuor pour la fin du Temps}), 
dont les rythmes sont   non r\'etrogradables,
et qui
sont respectivement 
$$3,\  5,\   8,\   5,\   3$$
$$4,\   3,\   7,\   3,\   4$$
$$ 2,\   2,\   3,\   5, \  3,\   2,\   2$$
$$1, \  1, \  3, \  2,\   2,\   1,\   2, \  2, \  3,\   1,\   1$$
$$2, \  1,\   1,\   1, \  3, \  1,\   1,\   1,\   2$$
$$ 2,\   1, \  1,\  1,\   3,\   1,\   1, \  1,\   2$$
$$ 1, \  1,\   1, \  1,\   1,\   3, \  1,\   1,\   1, \  1,\   1$$
$$3, \ 5, \ 8, \ 5, \ 3$$
(l'unit\'e de base \'etant la double croche).
Messiaen note dans son analyse que les valeurs totales des dur\'ees dans les
mesures
3
et 14 sont de 19 double croches, et dans les mesures 5, 6 et 7 cette
valeur est de 13 double croches, en faisant remarquer que 19 et 13 sont des
nombres premiers.

\medskip 

Il est naturel pour qui ne conna\^\i t pas l'\oe uvre  de Messiaen de se
demander pourquoi les rythmes non r\'etrogradables sont int\'eressants. Messiaen
r\'epond lui-m\^eme   \`a  cette question, et il donne deux r\'eponses, la premi\`ere   d'ordre esth\'etique
et la seconde, philosophique aussi, concernant le Temps. Dans sa \emph{Technique de mon langage musical}\footnote{Messiaen 1944.}, Messiaen parle du {\it charme} que peut produire un rythme non
r\'etrogradable sur l'auditeur de sa musique. En fait, il  parle en m\^eme
temps de ce charme et de celui produit par ses {\it modes \`a transpositions
limit\'ees}. M\^eme si l'on n'a pas encore rappel\'e la d\'efinition de ces modes, on peut signaler que Messiaen  d\'ecrit
les rythmes non r\'etrogradables et les 
modes \`a transpositions limit\'ees  comme des {\it impossibilit\'es
math\'ematiques}, impossibilt\'e  dans le domaine rythmique, car, comme il dit,
il est impossible de les r\'etrograder, puisqu'en les r\'etrogradant on obtient
le même rythme, et impossibilité dans le domaine modal, puisqu'on ne peut transposer au-delà d'un certain petit nombre de fois\footnote{« Un point fixera d'abord notre attention : {\it le charme des
impossibilit\'es}... Ce charme, \`a la fois voluptueux et contemplatif, r\'eside
particuli\`erement dans certaines impossibilit\'es math\'ematiques des domaines
modal et rythmique. Les modes, qui ne peuvent se transposer au-del\`a d'un
certain nombre de transpositions parce que l'on retombe toujours dans les
m\^emes notes ; les rythmes, qui ne peuvent se r\'etrograder parce que l'on
retrouve alors le m\^eme ordre des valeurs... » (Messiaen 1944, vol.1 p. 5).}.

\medskip

\`A la page 13 du m\^eme trait\'e, Messiaen d\'ecrit   les impressions
que font ces impossibilit\'es sur l'auditeur de sa musique :
   « Pensons maintenant \`a l'auditeur de notre musique modale et rythmique ; il
n'aura pas le temps, au concert, de v\'erifier les non-transpositions et les
non-r\'etrogradations, et \`a ce moment-l\`a, ces questions ne l'int\'eresseront
plus : \^etre s\'eduit, tel sera son unique d\'esir. Et c'est pr\'ecis\'ement ce
qui se produira : il subira malgr\'e lui le charme \'etrange deux impossibilit\'es : 
un certain effet d'ubiquit\'e tonale de la non-transposition, une certaine
unit\'e de mouvement (o\`u commencement et fin se confondent parce qu'ils sont
identiques) dans la non-r\'etrogradation, toutes choses  qui l'am\`eneront
certainement \`a cette sorte d'``arc-en-ciel th\'eologique" qu'essaie d'\^etre le
langage musical dont nous cherchons \'edification et th\'eorie. »

\medskip 

Une autre raison  pour laquelle   les rythmes
non r\'etrogradables étaient chers \`a Messiaen est r\'evélée 
dans le tome 2 de son \emph{Traité}. Qu'on le lise de droite \`a gauche ou de gauche \`a droite, dit-il, un tel rythme ne change pas, ce qui fait de lui un rythme  qui n'a ni commencement ni fin. Messiaen dit que l'une des {\it forces} de ces
rythmes réside dans le fait que « comme le Temps, le rythme non
r\'etrogradable est irr\'eversible. Il ne peut pas revenir en arri\`ere, sous
peine de se r\'ep\'eter... L'avenir et le pass\'e sont en miroir l'un par
rapport \`a l'autre ».

\medskip

 Pour terminer avec la notion de rythme, disons un mot sur le symbolisme du nombre dans le rythme. Bien que cela ne fasse pas partie des mathématiques proprement dites, il est bien connu que d'anciennes écoles mathématiques étaient attachées au symbolisme du nombre (on pense en particulier à celle de Pythagore). La rythmique hindoue foisonne de symbolisme du nombre. Par
exemple, un rythme \`a quatre temps, comme 1, 1, 1, 1 (T\^ala No. 99) s'appelle
{\it gaja}, qui en sanscrit veut dire   \'el\'ephant\footnote{Dans le
tableau qu'il a dress\'e des 120 {\it De\c c\^\i-T\^alas}, Messiaen indique
toujours le nom en sanscrit, la traduction fran\c caise de ce nom, puis les
caract\'eristiques rythmiques, symboliques et autres.}. Messiaen
explique\footnote{Messiaen 2001,  t. 1, p. 299.} que dans la culture hindoue, l'\'el\'ephant est la
manifestation de la force physique. « Ses quatre lourdes pattes et sa lourde et 
puissante d\'emarche sont repr\'esent\'ees par quatre dur\'ees ». L'\'el\'ephant
est symbolis\'e aussi par d'autres    rythmes \`a quatre dur\'ees, par
exemple 1, 1, 1, 3/2  (T\^ala No. 18), appel\'e {\it gajal\^\i la},
c'est-\`a-dire « jeu de l'\'el\'ephant », la derni\`ere valeur, plus longue que
les trois autres, reproduisant dans ce cas la lourdeur de la d\'emarche   de
l'\'el\'ephant\footnote{Messiaen 2001, t. 1, p. 276.}.
Comme autre exemple , signalons le T\^ala No.
105, compos\'e de 7 dur\'ees, successivement  2, 2, 2, 3, 3, 3, 1 et qui 
s'appelle {\it Candrakal\^a}, ce qui veut dire « beaut\'e de la lune ». Messiaen\footnote{Messiaen 2001, t. 1, p. 309.} pense que ce rythme comprend trois parties, d'abord  2, 2, 2
symbolisant la terre, puis 3, 3, 3 symbolisant le soleil, et enfin  1
symbolisant la lune.

\medskip 

Le T\^ala No. 27, que l'on a d\'ej\`a mentionn\'e,
1, 3, 2, 3, 3, 3, 2, 3, 1, 3, dans lequel s'entrecroisent une suite constante et une
suite croissante puis d\'ecroissante, a pour Messiaen  une
« importance extr\^eme »\footnote{Messiaen 2001, p. 317.}; il  repr\'esente,
selon lui, l'« union du Temps et de l'Eternit\'e ». D\'ej\`a utilis\'e
consciemment ou inconsciemment, par Stravinski, dans le {\it Sacre du
Printemps}, Messiaen dit que ce rythme est \`a la base de sa th\'eorie des
{\it personnages rythmiques}.

\medskip

 Il faut mentionner aussi \emph{le contrepoint de rythme}, une théorie que Messiaen utilisa dans ses pièces et qu'il développa dans ses trait\'es th\'eoriques. \`A l'origine, le  contrepoint est l'art de transformer et de
superposer des lignes m\'elodiques,  et Messiaen en a fait là aussi, une théorie propre aux rythmes.

\medskip

Classiquement, il y a trois op\'erations \'el\'ementaires dans le contrepoint mélodique : la rétrogradation,
l'inversion et la transposition. Parmi ces opérations, seule la première affecte le rythme.  En  contrepartie, il y a d'autres 
transformations de contrepoint qui sont propres au rythme, en particulier l'{\it augmentation} et la {\it diminution}. Elles  consistent à changer les valeurs des dur\'ees d'un certain motif mélodique en les multipliant
 par un facteur constant (on dit que l'on a une {\it imitation par
augmentation} si ce facteur est plus grand que 1, et {\it par diminution} s'il en est plus petit).
Avec Messiaen, l'augmentation et la diminution ont acquis un caract\`ere plus abstrait 
que dans le contrepoint classique parce qu'elles  sont appliqu\'ees
aux valeurs rythmiques seulement, c'est-à-dire sans que le nouveau motif mélodique ait forcément une relation avec l'ancien.
 Citons par exemple  la pi\`ece  No. V des {\it Vingt
Regards sur l'Enfant J\'esus} (1944) o\`u l'on trouve, au d\'ebut de la pi\`ece,  des
augmentations de rythmes par multiplication par  3/2, alors que le motif et son augmentation sont  compl\`etement diff\'erents du point
de vue m\'elodique ; ils sont même \'ecrits dans des modes diff\'erents, et ils forment ce
que Messiaen appelle un {\it canon rythmique}.

Plusieurs transformtions de rythmes pr\'eservant la propri\'et\'e d'\^etre
 non r\'etrogradables, sont d\'ecrites par Messiaen dans le
tome 2 du \emph{Traité} (p. 41). Ce sont celles d'{\it amplification sym\'etrique} et
d'{\it \'elimination sym\'etrique des extr\^emes}. L'amplification
sym\'etrique consiste \`a ajouter  \`a une formule rythmique donn\'ee une
autre formule et son renversement, respectivement d'un c\^ot\'e et de l'autre de
la formule initiale. Par exemple, au No. XX des {\it Vingt Regards sur l'Enfant
J\'esus}, le premier th\`eme, expos\'e \`a
la mesure 2, est un th\`eme tr\`es court dont la formule rythmique est 2, 1, 2 
(rythme non r\'etrogradable, o\`u l'unit\'e est ici aussi la double croche ; en
d'autres termes on a la suite croche, double croche, croche).  Ce th\`eme est
amplifi\'e  \`a la mesure 4, o\`u il devient 2, 2, 2, 1, 2, 2, 2. Il est
amplifi\'e de nouveau, de mani\`ere diff\'erente \`a la mesure 6, o\`u il devient
2, 2, 2, 1, 2, 1, 2, 2, 2.  On retrouve ce m\^eme rythme, avec les m\^emes
amplifications, plus loin dans la pi\`ece (mesure 82) : $$2 \ 1\  2$$ $$ 2\  2\ 
2\  1\  2\  2\  2$$ $$2\  3/2\  2\  2\  1\  2\  2\  3/2\  2$$
 
On peut mentionner aussi les op\'erations d'{\it agrandissement} et de {\it diminution de
la valeur centrale}, qui toutes deux pr\'eservent le caract\`ere non
r\'etrogradable. Plusieurs exemples sont donnés dans le tome 2 du \emph{Traité}.

\medskip 

Les techniques du contrepoint classique sont utilis\'ees 
en particulier pour \'ecrire des {\it canons}, et Messiaen  a inclus dans ses
pi\`eces  des {\it canons de rythmes}.  L\`a aussi, les exemples sont multiples
et plusieurs pi\`eces des {\it Vingt Regards} contiennent des {\it canons de
rythmes non r\'etrogradables}. Par exemple, la pi\`ece No. 5 en contient plusieurs.  Sur la partition de l'\oe uvre, Messiaen a
indiqué certains de façon explicite, pour en faciliter
l'analyse, les sujets et contre-sujets de ces canons, ainsi que les
diff\'erentes transformations qu'il leur a appliqu\'ees.

\medskip 

Retenons aussi que m\^eme si l'\'ecriture de Messiaen poss\`ede une rigueur tout à fait
math\'ematique, et  m\^eme si Messiaen demande \`a son lecteur et \`a son
interpr\`ete de « lire et ex\'ecuter exactement les valeurs marqu\'ees »\footnote{Messiaen 1944 chap. II.} sa musique exclut la monotonie. Dans ses dialogues avec Pierre Samuel\footnote{Samuel 1999, p. 102.}, il dit : « Une musique rythmique est une musique qui exclut la
r\'ep\'etition, la carrure et les divisions \'egales, qui s'inspire en somme
des mouvements de la nature, mouvements de dur\'ees libres et in\'egales », et
un peu plus loin (p. 103) : « La musique militaire est la n\'egation du rythme ».
Pour Messiaen\footnote{Messiaen 2001, p. 42.}, si dans le rythme il y a une p\'eriodicit\'e, c'est
« la vraie p\'eriodicit\'e, celle des vagues de la mer, le contraire d'une
r\'ep\'etition pure et simple. Chaque vague est diff\'erente de la
pr\'ec\'edente et de la suivante, par son volume, sa hauteur, sa dur\'ee, la
lenteur et la bri\'evet\'e de sa formation, la puissance de son climax, la
prolongation de sa chute, de son \'ecoulement, de son \'eparpillement ».

\medskip

Pour conclure cette section sur les rythmes, notons que m\^eme si Messiaen parle
d'un c\^ot\'e de rythmes qui sont « pens\'es pour le seul plaisir intellectuel du
nombre » (t. 1 p. 51), un plaisir qui peut sembler froid et d\'epourvu de
lyrisme, il compare d'un autre c\^ot\'e  la beaut\'e de certains rythmes qui
sont construits de mani\`ere rigoureuse, sur des structures math\'ematiques, \`a
la beaut\'e  de certains visages,   \`a celle
des jardins \`a la fran\c caise, \`a celles des  cath\'edrales romanes et à celle des ailes des papillons.

\section{ Les modes à transpositions limit\'ees}

  La seconde {\it impossibilit\'e math\'ematique} dont parle Messiaen dans le passage qu'on a cité de sa {\it Technique de mon langage musical} est celle des {\it
modes \`a transpositions limit\'ees}.

\medskip

Rappelons d'abord le sens du mot {\it
mode}. Messiaen, dans  
son \emph{Traité}, fait un exposé de la relation
entre mode et couleur. Adhérer complètement à son explication requiert  
une certaine facult\'e mentale (que Messiaen avait) : celle de pouvoir associer des couleurs à des notes ou ensemble de notes\footnote{voir, par exemple, Samuel 1999, p. 95, le chapitre intitul\'e {\it des sons et des couleurs}.}.
Faute de cela, on est   oblig\'e d'avoir recours \`a  une d\'efinition plus
terre-\`a terre, d'un mode comme une suite de notes musicales distinctes, qui
permet de d\'ecrire l' « atmosph\`ere » d'une  \oe uvre, ou
celle de l'endroit de l'\oe uvre o\`u l'on se trouve (une \oe uvre pouvant \^etre
polymodale). Pour utiliser un langage  « math\'ematique », un mode donne une
premi\`ere approximation de l'\oe uvre. Les {\it tonalit\'es} majeure  et mineure
sont des approximations  (assez rudimentaires) d'une telle notion, qui
s'appliquent \`a la musique {\it tonale} c'est-à-dire celle de l'occident d'après la Renaissance, 
 et jusqu'au d\'ebut du 20\`eme si\`ecle. La Gr\`ece  
et l'Inde antiques connaissaient des modes vari\'es et riches. Les
modes du plain-chant et ceux du chant byzantin trouvent
leurs racines dans ceux du chant grec antique.  Un retour à la musique modale, que l'on voit apparaître déjà vers la fin du dix-neuvième siècle dans les \oe uvres de Fauré, Satie et Debussy, met à la disposition du compositeur l'infinie variété de ces modes anciens.

\medskip

Il est naturel d'exprimer ces modes (et leurs transpositions) dans le langage mathématique du groupe des entiers modulo 12, de même qu'il est usuel d'expliquer les opérations de transposition en utilisant ce même langage. 
\medskip

Techniquement parlant, un {\it mode \`a transpositions limit\'ees} est un ensemble de notes de la gamme tempérée à 12 sons
ayant  la propri\'et\'e que quand on transpose cette suite par
un certain intervalle « petit » en un sens à préciser, on retombe sur la m\^eme suite, les notes \'etant consid\'er\'ees \`a octave pr\`es\footnote{Rappelons Rameau est le premier à avoir pris dans ses axiomes de la musique et théorisé l'identification de deux notes à une octave près (Rameau 1750), même si cette  identification existait, du point de vue de la pratique, longtemps avant lui (les hommes et les femmes chantant simultanément « à l'unisson » alors qu'ils chantent à une octave près).}. 

\medskip On commence par représenter les 12 notes de la gamme chromatique, {\it
do, do dièse, r\'e, r\'e dièse, mi, fa, fa dièse, sol, sol dièse, la, la dièse, si}, par la suite des  nombres entiers (modulo 12) 0, 1, 2, 3, ...
11. 
Un mode, étant une
suite de notes distinctes parmi les 12 notes de la gamme chromatique définies à une octave près, devient un sous-ensemble de l'ensemble
$0,\ 1,\ ...,\ 11$. 

L'ensemble des entiers modulo 12 est muni d'une opération d'addition naturelle, le résultat de l'addition de deux nombres étant soit le nombre obtenu par l'opération d'addition ordinaire, ou bien (si ce nombre dépasse 11), ce même nombre auquel on a retranché 12.

L'opération de {\it transposition} au niveau des notes, correspond alors, au
niveau de la représentation des notes par des entiers modulo 12, \`a une translation dans le groupe des entiers modulo 12 (appelée translation modulo 12).

\medskip 

Un
{\it mode \`a transpositions limit\'ees à sa premi\`ere transposition} devient doc un sous-ensemble cycliquement ordonnée de l'ensemble $ 0,\ 1,\ ...,\ 11 $ commençant par 0 et invariante par une « petite » translation modulo 12, à permutation circulaire près.  Messiaen retient sept modes à transpositions limitées particuliers, que l'on va d\'ecrire
maintenant, \`a leur premi\`ere transposition. Notons tout de suite que cette
classification en sept modes est distincte de toutes celles des  modes connus avant
Messiaen (que ce soit les modes de la Gr\`ece antique, de l'Inde, de la Chine ou
du chant gr\'egorien).  

La liste des sept modes de Massian est la suivante :

\medskip 
\noindent {\it Premier mode.---} Il est « deux
fois transposable », c'est-\`a-dire, en langage math\'ematique,  invariant par la
translation par 2. Du point de vue musical, c'est la {\it gamme de tons}, qui fut largement utilis\'ee par Debussy (et que l'on appelle parfois la
« gamme de tons de Debussy »). En d'autres termes, c'est la suite de
notes  $$\hbox{\it do, r\'e, mi,  fa dièse, sol dièse, la dièse}.$$ 
Utilisant le langage
 des entiers modulo 12, c'est la suite 
$$0,\  2,\ 4,\ 6,\ 8,\ 10.$$
\`A une permutation circulaire pr\`es, elle est invariante par la translation par 2 modulo 12.

\medskip

\noindent {\it Second  mode.---} Ce mode est « trois fois transposable », 
 c'est-\`a-dire invariant par 
translation par 3. Il est repr\'esent\'e par la suite de
notes 
$$ \hbox{\it do,  do dièse, ré dièse,   mi,  fa dièse,   sol, 
la, la dièse}$$ ou, de mani\`ere \'equivalente, par la suite 
d'entiers modulo 12 $$ 0,\ 1,\ 3, \ 4,\ 6,\ 7,\ 9,\  10.$$
Messiaen utilisait ce mode de façon systématique. Il dit qu'on le trouve « \`a l'\'etat de timide \'ebauche » dans la musique de  Rimski-Korsakov, Scriabine,
Ravel et Stravinski.  

\medskip

\medskip 

\noindent {\it Troisi\`eme  mode.---} Ce mode est « quatre fois transposable », 
 c'est-\`a-dire invariant par la
translation par 4.
 Il est repr\'esent\'e par la suite de
notes 
$$\hbox{\it do, r\'e, ré dièse, mi, fa dièse, sol, sol dièse,
la dièse, si} $$ ce qui correspond à la suite
d'entiers modulo 12 $$0,\ 2,\ 3,\ 4,\ 6,\ 7,\ 8,\ 10,\ 11 .$$

\medskip 
\medskip
 Les quatre derniers modes (les modes 4, 5, 6 et 7) sont « six fois
transposables », 
 c'est-\`a-dire que chacun d'eux est invariant par la
translation  par 6 modulo 12. Ils sont repr\'esent\'es respectivement
par les notes, ou suites d'entiers suivants :

\medskip 
\medskip 

\noindent {\it Quatri\`eme  mode.---}
$$\hbox{\it do, do dièse, r\'e, fa, fa dièse, sol, sol dièse, si 
}$$
 ou $$0, \ 1, \ 2, \ 5, \ 6,\ 7,\ 8,\ 11 .$$

\noindent {\it Cinqui\`eme  mode.---}
$$\hbox{\it do, do dièse, fa, fa dièse, sol, si}$$
ou
$$0,\ 1,\ 5,\ 6,\ 7,\ 11 .$$

\noindent {\it Sixi\`eme   mode.---}

$$\hbox{\it do, r\'e, mi, fa, fa dièse, sol dièse, la dièse, si}$$
ou
$$0,\ 2,\ 4,\ 5,\ 6,\ 8,\ 10,\ 11 .$$

\noindent {\it Septi\`eme   mode.---}

$$\hbox{\it  do, do dièse, r\'e, ré dièse, mi, fa dièse, sol,
sol dièse, la, si}$$ ou
$$0,\ 1,\ 2,\ 3,\ 4,\ 6,\ 7,\ 8,\ 9,\ 11 .$$

\medskip

La liste des modes
\`a transpositions limit\'ees s'arr\^ete aux modes six fois
transposables car 6 est le plus grand diviseur de 12 qui ne soit pas égal à 12. 

Cette liste de 7
modes n'est pas exhaustive par rapport aux propri\'et\'es math\'ematiques 
que l'on a mentionn\'ees. Messiaen d'ailleurs donne lui-même une liste de suites de
notes de la gamme chromatique qui ne sont pas des transpositions des modes 4, 5,
6 ou 7, et qui sont pourtant six fois transposables\footnote{Messiaen 2001, p. 54.}, mais il les considère comme des {\it modes tronqu\'es}, et non pas de vrais modes. 

Dans la correspondance mode-couleur établie par Messiaen, chacun des
sept modes qu'il  considère, 
 \`a chacune  de ses transpositions, correspond 
\`a une certaine couleur bien définie. Nous avons mentionné le « charme des impossibilit\'es » des modes à transpositions limitées. Selon Messiaen, ces modes « r\'ealisent dans le sens vertical (transposition) ce
que les rythmes non r\'etrogradables r\'ealisent dans le sens horizontal
(r\'etrogradation) »\footnote{Messiaen 1944, p. 13.}.  Math\'ematiquement parlant, toutes ces impossibilit\'es
s'expriment par des sym\'etries. La notion d'nvariance par un groupe de transformations est équivalent à celle de symétrie.

Il existe une troisi\`eme « impossibilit\'e mathématique »
dont   Messiaen ne parle pas dans sa \emph{Technique}, mais il le fait quelques ann\'ees
plus tard, dans son \emph{Traité} :  c'est celle des {\it permutations sym\'etriques}, dont on va parler maintenant.

 \section{Permutations sym\'etriques}

Messiaen appelle {\it permutation sym\'etrique} une permutation (d'un nombre
fini d'objets) ayant la propriété que quand on la répète un petit nombre de fois, on retombe sur l'identité. En termes mathématiques, l'ordre du groupe cyclique que cette permutation engendre est petit. C'est la même notion qui intervient dans les modes à transpositions limitées. Le sens de l'adjectif « petit » peut être pris dans le sens de pouvoir être décelé par un auditeur de la musique.

Les
permutations, en général, jouent un r\^ole déterminant dans la musique de Messiaen. \`A partir
d'une suite de notes, de dur\'ees, de nuances, etc., on peut en obtenir une
autre par permutation; par exemple, à partir de la suite de durées \emph{croche noire noire}, on obtient, par  permutation circulaire vers la droite, \emph{noire noire croche} ; en appliquant de nouveau la même permutation, on obtient \emph{noire croche  noire}, et une troisième fois la même, \emph{croche noire noire}. Ainsi, la permutation que l'on itère ici est d'ordre trois.  La  r\'etrogradation (le fait de lire la suite en sens inverse) est un cas particulier de
permutation ; elle est d'ordre deux. Une permutation peut s'appliquer au rythme comme à d'autres paramètres de la musique. Pour que l'effet de \emph{toutes} les permutations d'une suite finie d'objets soit exploitable en musique (ou, du moins, pour qu'il soit perceptible à l'oreille), il faut que le nombres d'objets  soit assez petit -- ne dépassant  pas 4 ou 5. D\`es que le nombre d'objets est grand, le nombre total de permutations est trop grand pour que cet
ensemble soit utilisable.  C'est pour cela que Messiaen propose de se
limiter à un petit nombre de permutations, en utilisant des permutations d'ordre petit. C'est ce qu'il entend par « permutation symétrique ».

\medskip

Dans ce cadre, une pièce que l'on donne souvent en exemple est la {\it Chronochromie} (Couleur du temps), dans laquelle Messiaen applique une permutation symétrique à un  rythme à 32 dur\'ees. Il part du rythme d\'efini par la suite   
chromatique de dur\'ees, commençant  à la triple croche et not\'ee 1 ; 2 ; 3 ; 4 ; ... ; 32. 
(La première valeur est la triple croche, la seconde est 2 triples croches, etc.)
Messiaen utilise la permutation suivante :
3 ; 28 ; 5 ; 30 ; 7 ; 32 ; 26 ; 2 ; 25 ; 1 ; 8 ; 24 ; 9 ; 23 ; 16 ; 17 ; 18 ;
22 ; 21 ; 19 ; 20 ; 4 ; 31 ; 6 ; 29 ; 10 ; 27 ; 11 ; 15 ; 14 ; 12 ; 13.
Ainsi, appliquée à notre suite chromatique de durées, cette permutation donne une suite de dur\'ees valant, dans l'ordre, 
3 triples croches, 28 triple croches, 5 triple croches, etc.
La m\^eme permutation est appliquée à la suite obtenue, et ainsi de suite.
Dans tome 3 de son \emph{Traité}, p. 15 \`a 66, Messiaen dresse un tableau complet de la suite des
itérés ainsi obtenue ; il en ressort que cette suite n'a que 36 \'el\'ements.
En d'autres termes, en appliquant la m\^eme permutation 36 fois, on retombe sur
la suite initiale de durées : 1 ; 2 ; 3 ; ... ; 32.

\medskip

Il existe un type de permutations auquel
Messiaen attache une  importance particuli\`ere  (et qui est peut-être à l'origine de l'expression « permutation sym\'etrique ») : on part du milieu d'une suite d'objets (par exemple d'un rythme) et l'on prendre ensuite dans l'ordre  un objet de part et d'autre, en partant du centre jusqu'aux extrémités. Une telle permutation est d'ordre petit. 

Par exemple , pour
une permutation de trois objets not\'es 1 ; 2 ; 3, en faisant la mouvement
dans le sens gauche-droite, on obtient la permutation 
2 ; 1 ; 3. En itérant une fois, on retombe sur 1 ; 2 ; 3.
La permutation considérée est donc d'ordre deux.

Pour un ensemble de 4 objets, ce procédé de permutation sym\'etrique donne successivement les suites
 (1  ; 2 ; 3  ;  4), (2  ; 3 ; 1 ;  4),  (3  ; 1 ; 2  ;  4),  (1  ; 2 ; 3  ;  4).
La permutation est donc d'ordre 4.

\medskip

 Messiaen résume avec les mots suivants ce procédé d'obtention de permutations d'ordre petit :  « Il s'agit d'un proc\'ed\'e qui  
correspond exactement \`a ce que j'appelle, dans les modes, les modes \`a
transpositions limit\'ees, et, dans les rythmes, les rythmes non
r\'etrogradables. C'est un proc\'ed\'e qui repose sur une impossibilit\'e. Ce
sont des dur\'ees qui se succ\`edent dans un certain ordre et que l'on relit
toujours dans l'ordre de d\'epart; prenons, par exemple, une gamme chromatique
de trente-deux dur\'ees : on les intervertit selon un ordre choisi, on
num\'erote le r\'esultat de un \`a trente-deux, puis on relit ce r\'esultat dans
le premier  ordre, et ainsi de suite jusqu'\`a ce  qu'on retrouve la gamme
chromatique des trente-deux dur\'ees de d\'epart. Ce syst\`eme produit  des
rythmes int\'eressants et tr\`es \'etranges, mais il pr\'esente surtout
l'avantage d'\'eviter  un nombre absolument fabuleux de permutations. Vous savez
qu'avec le chiffre douze, tellement aim\'e des s\'eriels, le nombre des
permutations est 479 001 600 ! Il faudrait des ann\'ees pour les \'ecrire.
Tandis qu'avec mon proc\'ed\'e, on peut, avec des chiffres plus importants --
trente-deux ou soixante-quatre -- obtenir les meilleures permutations, supprimer
les permutations secondaires qui n'aboutissent qu'\`a des r\'ep\'etitions, et
travailler sur un nombre de permutations raisonnable, pas tr\`es loin du chiffre
de d\'epart »\footnote{Samuel 1999, p. 119.}. Messiaen explique ainsi le nom de la pièce  : « Les dur\'ees et les permutations de dur\'ees rendues sensibles par des colorations sonores : c'est bien une « couleur du temps, une chronochromie »\footnote{Samuel 1999, p. 222.}.
 
 \medskip
 
 On pourrait parler aussi des mathématiques sous-jacentes aux techniques que l'on a appelées « super-sérielles » dans la musique de Messiaen, qui, généralisent les techniques sérielles développées au départ par Arnold Schoenberg et son école, pour les séries à 12 sons. Messiaen utilisa ces techniques de façon épisodique et il ne les  considéra pas comme importantes. Il déclara notamment : « J'ai \'et\'e tr\`es contrari\'e de l'importance
absolument d\'emesur\'ee que l'on a accord\'ee \`a une petite \oe uvre, qui n'a
que trois pages et qui s'appelle {\it Mode de valeurs et d'intensit\'es}, sous le
pr\'etexte qu'elle aurait \'et\'e \`a l'origine de l'\'eclatement s\'eriel dans
le domaine des attaques, des dur\'ees, des intensit\'es, des timbres, bref, de
tous les param\`etres musicaux. Cette musique a peut-\^etre \'et\'e
proph\'etique, historiquement importante, mais musicalement, c'est trois fois
rien... »\footnote{Samuel 1999,  p. 72.}.

 \medskip

On peut parler de nouveau de symbolisme du nombre.

\medskip

 Indépendamment des techniques mathématiques qu'il a introduites, Messiaen aimait les nombres et leur attachait des valeurs symboliques, dans la tradition des Pythagoriciens, qui sont les fondateurs de la première école mathématique connue et digne de ce nom.  Le tome 3 de son \emph{Traité} contient plusieurs pages concernant ce sujet, dans lesquelles il explique par exemple que le nombre 1 symbolie l'unité divine,  le nombre 3 la Sainte Trinité, le nombre 7 la création, et le nombre 10 la décade pythagoricienne, qui en même temps est Dieu et le monde (3+7). Comme le fit Jean-Sébastien Bach avant lui, Messiaen a introduit ces symboles dans sa musique. Il nous dit par exemple :  « La {\it Transfiguration} se divise en deux groupes de sept
pi\`eces : j'ai adopt\'e le chiffre symbolique de sept, que l'on trouve
constamment dans la Bible, et en particulier dans l'Apocalypse »\footnote{Samuel 1999, p. 232.}.
En parlant des {\it Vingt Regards sur l'Enfant J\'esus} (texte accompagnant
l'enregistrement  Erato, avec Yvonne Loriod au piano), il écrit : « Les num\'eros des
pi\`eces sont ordonn\'es par les contrastes de temps, d'intensit\'e, de
couleur -- et aussi par des raisons symboliques... » Par exemple, le  {\it Regard du Temps} porte le num\'ero 9, car « le Temps a vu
na\^\i tre en lui celui qui est \'eternel, en l'enfermant dans les 9 mois de la
maternit\'e », etc. Dans la pr\'eface de la partition du {\it Quatuor pour la fin du Temps} (Ed.
Durand), Messiaen explique le choix de la num\'erotation et du nombre total
des mouvements de cette pièce (8) par la symbolique de ces nombres.

\section{Conclusion}

S'il passait beaucoup de temps à compter, transformer et réfléchir sur des suites de nombres. S'il utilisait consciemment les mathématiques dans ses compositions, ce n'est pas seulement dans le but de produire des effets nouveaux. Il aimait cette science parce qu'il sentait qu'elle fait  partie d'un monde abstrait immuable, en dehors du Temps. C'est pour cette raison qu'il accorde au nombre une place importante dans son  \emph{Traité}, au même titre que la m\'emoire, l'\'eternit\'e et le cosmos.  
 
   Artiste, observateur de la nature -- les oiseaux en particulier -- et aimant les mathématiques, Messiaen nous fait penser à Léonard de Vinci, qui, à un moment de sa vie a délaissé la peinture pour approfondir ses connaissances en mathématiques. Gabriel Séailles, son biographe connu, dans son livre \emph{Léonard de Vinci, l'artiste et le savant : 1452-1519 : essai de biographie psychologique}, cite une lettre du Révérend  Petrus de Nuvolaria à la Duchesse de Milan Isabelle d'Esté,  une figure connue de la Renaissance italienne, dans laquelle il écrit, à propos de l'artiste : « Ses études mathématiques l'ont à ce point dégoûté de la peinture, qu'il supporte à peine de prendre une brosse. »  Séailles cite aussi  Sabba da Castiglione, un auteur humaniste contemporain de Leonardo, qui écrit à propos de ce dernier:  « Quand il devait se consacrer à la peinture, où sans aucun doute il eût été un nouvel Apelles, il se donna tout entier à la géométrie, à l'architecture, à l'anatomie »\footnote{Séailles 1892.}. 
   
   \medskip
   
 Le \emph{Trait\'e de Rythme, de couleur et d'ornithologie}  de Messiaen nous fait penser aux écrits d'Aristote par leur leurs classifications exhaustives et le sens inouï du détail, avec leur analyse fine de la perception humaine et les réflexions qu'ils contiennent sur le temps, la durée, la nature, l'éternité, sur la différence entre ce qui est immuable et indivisible et ce qui est changeant et ce qui est périodique « par alternances du même et de l'autre». Cette proximité avec Aristote est encore un trait commun à Messiaen et Thom. Ce dernier fut en effet, durant les dix ou quinze dernières années de sa carrière, le grand défenseur des idées du Stagirite\footnote{Voir sur le sujet l'exposé dans Papadopoulos 2019A.}.

\`A propos de l'aspect cosmique de la pensée de Messiaen, en relation avec les mathématiques auxquelles il attachait tant d'importance, on peut aussi citer le passage bien connu de l'« Essayeur » de Galilée, dans lequel ce dernier considère que le « livre de l'univers »  est écrit en termes mathématiques : « La philosophie est écrite dans cet immense livre qui continuellement reste ouvert devant les yeux (ce livre qui est l'Univers), mais on ne peut le comprendre si, d'abord, on ne s'exerce pas à en connaître la langue et les caractères dans lesquels il est écrit. II est écrit dans une langue mathématique, et les caractères en sont les triangles, les cercles, et d'autres figures géométriques, sans lesquelles il est impossible humainement d'en saisir le moindre mot; sans ces moyens, on risque de s'égarer dans un labyrinthe obscur »\footnote{trad. Chauviré 1980, p. 141.}.

\medskip

L'influence de Messiaen sur l'utilisation des mathématiques en musique fur très grande.
Pendant sa longue carrière de pédagogue, il eut des élèves qui  qui, chacun à sa manière, utilisa les mathématiques dans ses compositions (on pense à Stockhausen, Boulez, Xenakis, etc.). 
 
On peut citer ici Messiaen, rencontrant le jeune Xenakis, qui lui explique qu'il voudrait suivre ses cours mais qu'il n'avait pas le bagage musical classique :  «  [...] c'était un homme
tellement hors du commun ! Je lui ai dit : `Non. Vous avez déjà trente ans, vous avez la chance d'être Grec,
d'avoir fait des mathématiques, d'avoir fait de l'architecture. Profitez de ces choses-là, et faites-les dans votre
musique'. Je crois finalement que c'est ce qu'il a fait »\footnote{Cité dans Matossian 1981, p. 58.}.

 Il est juste de conclure cet article en se posant la question de savoir ce que Messiaen considérait comme le plus important dans sa musique. La question lui fut posée par Claude Samuel, dans son entretien réalisé en 1967 c'est-à-dire en plein milieu de sa carrière de compositeur :  « Quelles `expressions' voulez-vous donc défendre en écrivant de la musique ? Quelles impressions voulez-vous communiquer à vos éditeurs ? » Messiaen répond : « La première idée que j'ai voulu exprimer, celle qui est la plus importante parce qu'elle est placée au-dessus de tout, c'est l'existence des vérités de la foi catholique.  J'ai la chance d'être catholique ; je suis né croyant et il se trouve que les textes sacrés m'ont frappé dès mon enfance. Un certain nombre de mes \oe uvres sont donc destinées à mettre en lumière les vérités théologiques de la foi catholique.C'est là le premier aspect de mon \oe uvre, le plus noble, sans doute le plus utile, le plus valable, le seul peut-être que je ne regretterai pas à l'heure de ma mort. Mais je suis un être humain, comme tous les êtres humains je suis sensible à l'amour humain que j'ai voulu exprimer dans trois de mes \oe uvres par l'intermédiaire du plus grand mythe de l'amour humain, Tristan et Iseult. Enfin, j'admire profondément la nature. Je pense que la nature nous surpasse infiniment et je lui ai toujours demandé des le\c cons ; par goût, j'ai aimé les oiseaux, j'ai donc spécialement interrogé les chants des oiseaux : j'ai fait de l'ornithologie. Il y a dans ma musique cette juxtaposition de la foi catholique, du mythe de Tristan et Yseult et l'utilisation excessivement poussée des chants d'oiseaux »\footnote{Samuel 1967, p. 11.}.
La musique de Messiaen est pleine de références théologiques, mais le côté religieux et méditatif de Messiaen sont probablement  trop bien connus pour que cela vaille la peine d'y insister ici.

\bigskip
\centerline{\bf Références bibliographiques}

\medskip

 \noindent Chauviré 1980 \\
Chauviré, Christiane, L'Essayeur de Galilée, Coll. Annales littéraires de l'Université de Besançon, no. 234, Paris, Les Belles Lettres, 1980.

\medskip

 \noindent Ide 1999 \\
Ide, Pascal, Olivier Messiaen, un musicien ébloui par l'infinité de Dieu, Nouvelle Revue Théologique, 121 (juillet-septembre 1999) n° 3, p. 436-453.

\medskip

 \noindent Mathon 1995 \\
Mathon,   Geneviève, Les couleurs de la cité Céleste, éléments d'analyse,  Musurgia
Vol. 2, No. 1 (1995), p. 141-158.

\medskip

 \noindent Matossian 1981 \\
Matossian, Nouritza, Iannis Xenakis, Paris, Fayard, 1981.

\medskip

 \noindent Messiaen 1944 \\
 Messiaen, Olivier, Technique de mon langage musical, ed. Alphonse Leduc,
Paris 1944.

\medskip

 \noindent Messiaen 1966 \\
Messiaen, Olivier,  \emph{Couleurs de la cité céleste}, Leduc, 1966.

\medskip

 \noindent Messiaen 2001 \\
 Messiaen, Olivier, Trait\'e de rythme, de couleur et
d'ornithologie, 7 tomes, ed. Alphonse Leduc, Paris.   

\medskip

 \noindent  Papadopoulos 2019A \\
  Papadopoulos, Athanase, Thom et Aristote : Penser la forme et la nature, In :  Ren\'e Thom, Portait mathématique et philosophique,   CNRS \'Editions, Paris, 2018, p.  203--260.

\medskip

 \noindent  Papadopoulos 2019B \\
Papadopoulos, Athanase, Logos et analogie : Notes sur deux thèmes récurrents dans la pensée de Thom In :  Ren\'e Thom, Portait mathématique et philosophique,   CNRS \'Editions, Paris, 2018, p. 261--292.
\medskip

 \noindent Rameau 1750 \\
 Rameau, Jean-Philippe, Démonstration du principe de l'harmonie, servant de base à tout l'art musical théorique et pratique, Durand et Pissot, Paris,  1750.

\medskip

 \noindent  Samuel 1967 \\
 Samuel, Claude, Entretiens avec Olivier Messiaen, Belfond, 1967.

\medskip

 \noindent Samuel 1999 \\
Samuel, Claude, Permanences d'Olivier Messiaen 
(Dialogues et commentaires), ed. Actes Sud 1999.

\medskip

 \noindent Séailles 1892 \\
Séailles, Gabriel,  \emph{Léonard de Vinci, l'artiste et le savant : 1452-1519 : essai de biographie psychologique}, Perrin, Paris, 1892.

\medskip

 \noindent {\sc Thom} 1988\\
 Thom, René, \emph{Esquisse d'une s\'emiophysique : Physique aristot\'elicienne et th\'eorie des catastrophes}, Paris,  InterEditions,  1988.

\end{document}